\documentclass[a4paper,12pt]{article}

\usepackage{amssymb}
\usepackage{amsmath}
\usepackage{mathptmx}
\usepackage{color}
\usepackage{cancel}
\usepackage{bbm}
\usepackage{pstricks}
\usepackage{psfrag}
\def\ap{'\thinspace}

\def\spacingset#1{\def\baselinestretch{#1}\small\normalsize}
\spacingset{1.3}

\newtheorem{lemma}{Lemma}[section]
\newtheorem{theorem}{Theorem}[section]
\newtheorem{remark}{Remark}[section]

\newtheorem{proposition}{Proposition}[section]

\def\be{\begin{equation}}
\def\ee{\end{equation}}
\def\bea{\begin{eqnarray}}
\def\eea{\end{eqnarray}}
\def\beann{\begin{eqnarray*}}
\def\eeann{\end{eqnarray*}}
\def\bsea{\begin{subeqnarray}}
\def\esea{\end{subeqnarray}}
\def\bmat{\left[ \begin{array}}
\def\emat{\end{array} \right]} 

\def\ap{'\thinspace}

\def\proof{\noindent{\bf{\em Proof:}\ \ }}
\def\QED{\mbox{\rule[0pt]{1.5ex}{1.5ex}}}
\def\endproof{\hspace*{\fill}~\QED\par\endtrivlist\unskip}
\newcommand{\real}{{\mathbb{R}}}
\newcommand{\complex}{{\mathbb{C}}}

\def\gD{{\cal D}}

\def\gR{{\cal R}}
\def\gS{{\cal S}}

\def\gU{{\cal U}}

\newcommand{\ima}{\operatorname{im}}
\newcommand{\rank}{\operatorname{rank}}

\newcommand{\defi}{\stackrel{\text{\tiny def}}{=}}
\def\tra{{\scalebox{.6}{\mbox{T}}}}
\def\bsmat{\left[ \begin{smallmatrix}}
\def\esmat{\end{smallmatrix} \right]}

\definecolor{Royalblue}{cmyk}{1,0.30,0.2,0.2}

\newcommand{\tp}{{^\tra}}

\begin{document}
\begin{titlepage}
\title{\vspace{-10mm}
A reduction technique for\\ Generalised Riccati Difference Equations\thanks{Partially supported by the Italian Ministry for Education and Research (MIUR) under PRIN grant n. 20085FFJ2Z ``New Algorithms and Applications of System Identification and Adaptive Control" and by the Australian Research Council under the grant FT120100604. Research carried out while the first author was visiting Curtin University, Perth (WA), Australia. }\vspace{10mm}}
\author{{\large Augusto Ferrante$^\ddagger$  \quad Lorenzo Ntogramatzidis$^\star$ }\\
       {\small $^\ddagger$Dipartimento di Ingegneria dell\ap Informazione}\\[-3mm]
        {\small    Universit\`a di Padova, via Gradenigo, 6/B -- 35131 Padova, Italy}\\[-3mm]
       {\small     {\tt augusto@dei.unipd.it}} \\ 
       {\small     $^\star$Department of Mathematics and Statistics}\\[-3mm]
       {\small     Curtin University, Perth WA, Australia.}\\[-3mm]
       {\small    {\tt L.Ntogramatzidis@curtin.edu.au}}
 }%
\thispagestyle{empty} \maketitle \thispagestyle{empty}

\begin{center}
\begin{minipage}{14.8cm}
\begin{center}
\bf Abstract
\end{center}
This paper proposes a reduction technique for the generalised Riccati difference equation arising in optimal control and optimal filtering.
This technique relies on a study on the generalised discrete algebraic Riccati equation. 
In particular, an analysis on the eigenstructure of the corresponding extended symplectic pencil enables to identify a subspace in which all the solutions of the generalised discrete algebraic Riccati equation are coincident. 
This subspace is the key to derive a decomposition technique for the generalised Riccati difference equation that isolates its nilpotent part, which becomes constant in a number of steps equal to the nilpotency index of the closed-loop, from another part that can be computed by iterating a reduced-order generalised Riccati difference equation.

\end{minipage}
\end{center}
\begin{center}
\begin{minipage}{14.2cm}
\vspace{2mm}
\end{minipage}
\end{center}
\thispagestyle{empty}
\end{titlepage}
%

\section{Introduction}
\label{secintro}
Consider the classic finite-horizon Linear Quadratic (LQ) optimal control problem. In particular, consider the discrete linear time-invariant system governed by the difference equation
\bea
\label{eqsys}
x_{t+1} = A\,x_t+B\,u_t,
\eea
where $A \in \real^{n \times n}$ and $B \in \real^{n \times m}$, and where, for all $t \ge 0$, $x_t\in \real^n$ represents the state and $u_t \in \real^m$ represents the control input. Let the initial state $x_0\in \real^n$ be given. The problem is to find a sequence of inputs $u_t$, with $t = 0,1, \ldots,T-1$, minimising the cost function
\bea
\label{cost}
J(x_0,u) \defi \sum_{t=0}^{T-1} \bmat{cc} x_t^\tra \;&\; u_t^\tra \emat \bmat{cc} Q & S \\ S^\tra & R \emat \bmat{c} x_t \\ u_t \emat+x_T^\tra\,P\,x_T.
\eea
We assume that the weight matrices $Q\in \real^{n \times n}$, $S \in \real^{n \times m}$ and $R \in \real^{m \times m}$ are such that the {\em Popov matrix} $\Pi$ is symmetric and positive semidefinite, i.e.,
\bea
\Pi \defi \bmat{cc} Q & S \\ S^\tra & R \emat =\Pi^\tra \ge 0.
\eea
We also assume that $P=P^\tra \ge 0$.
The set of matrices $\Sigma=(A,B,\Pi)$ is often referred to as {\em Popov triple}, see e.g. \cite{Ionescu-OW-99}. 
We recall that, for any time $t$, the set $\gU_t$ of all optimal inputs can be {parameterised in terms of an arbitrary $m$-dimensional signal $v_t$ as}
$\gU_t=\{-K_t\,x_t+G_t\,v_t\}$, where\footnote{The symbol $M^\dagger$ denotes the Moore-Penrose pseudo-inverse of matrix $M$.}
\bea
K_t & = & (R+B^\tra\,X_{t+1}\,B)^\dagger (S^\tra+B^\tra\,X_{t+1}\,A), \\
G_t & = & I_m-(R+B^\tra\,X_{t+1}\,B)^\dagger (R+B^\tra\,X_{t+1}\,B),
\eea 
in which $X_t$ is the solution of the Generalised Riccati Difference Equation GRDE($\Sigma$)
\bea
\label{grde}
X_{t}  =   A^\tra\,X_{t+1}\,A-(A^\tra\,X_{t+1}\,B+S)(R+B^\tra\,X_{t+1}\,B)^\dagger(B^\tra\,X_{t+1}\,A+S^\tra)+Q \qquad
\eea
 iterated backwards from $t=T-1$ to $t=0$ using the terminal condition
 \bea
 \label{term}
 X_T=P,
 \eea
 see \cite{Rappaport-S-71}.  The equation characterising the set of optimal state trajectories is 
\beann
x_{t+1}=(A-B\,K_t)\,x_t-B\,G_t\,v_t.
\eeann
The optimal cost is $J^\ast=x_0^\tra\,X_0\,x_0$.\\[-2mm]

Despite the fact that it has been known for several decades that the generalised discrete Riccati difference equation provides the solution of the classic finite-horizon LQ problem, 
this equation has not been studied with the same attention and thoroughness that has undergone the study of the standard discrete Riccati difference equation.
 The purpose of this paper is to attempt to start filling this gap. In particular, we want to show a reduction technique for this equation that allows to compute its solution by solving a smaller equation with the same recursive structure, with obvious computational advantages.
In order to carry out this task, several ancillary results on the corresponding generalised Riccati equation are established, which constitute an extension of those valid for standard discrete algebraic Riccati equations presented in \cite{Ferrante-W-07} and \cite{Ferrante-04}.
In particular, these results show that the nilpotent part of the closed-loop matrix is independent of the particular solution of the generalised algebraic Riccati equation. Moreover, we provide a necessary and sufficient condition expressed in sole terms of the problem data for the existence of this nilpotent part of the closed-loop matrix. This condition, which appears to be straightforward for the standard algebraic Riccati equation, becomes more involved -- and interesting -- for the case of the generalised Riccati equation.
{We then show that every solution of the generalised algebraic Riccati equation coincide along the largest  eigenspace associated with the eigenvalue at the origin of the closed-loop, and that this subspace can be employed to decompose the generalised Riccati difference equation into a nilpotent part, whose solution converges to the zero matrix in a finite number of steps {(not greater than $n$)} and a part which corresponds to a non-singular closed-loop matrix, and is therefore easy to handle with the standard tools of linear-quadratic optimal control.}
{ As a consequence, our analysis permits a generalisation of a long series of results 
aiming to the closed form representation of the optimal control, see \cite{Ferrante-N-05,Ferrante-N-06,EZ,Ferrante-Ntogramatzidis-IEEE-2012} and, for the continuous-time counterpart, \cite{Ferrante-M-N,Ferrante-Ntog-EJC-07,Ferrante-N-10}.
}
{
Our analysis of the GRDE is based on the general theory on generalised {\em algebraic} Riccati equation presented in \cite{Stoorvogel-S-98} and on some recent developments derived in \cite{Ferrante-N-12-sub,Ferrante-N-12-sub2}. } 
  
\section{The Generalised Discrete Algebraic Riccati Equation}
We begin this section by recalling two standard linear algebra results that are used in the derivations throughout the paper.
  {
  \begin{lemma}
  \label{lem1}
  Consider $P=\left[ \begin{smallmatrix} P_{11} & P_{12} \\[1mm] P_{12}^\tra & P_{22} \end{smallmatrix} \right]=P^\tra \ge0$. Then,
    \begin{enumerate}
  \item $\,\ker P_{12} \supseteq \ker P_{22}$; 
  \item $\,P_{12}\,P_{22}^\dagger\,P_{22}= P_{12}$;  
  \item $\,P_{12}\,(I-P_{22}^\dagger P_{22})=0$; 
  \item $\,P_{11}-P_{12} P_{22}^\dagger P_{12}^\tra \ge 0$; 
  \end{enumerate}
    \end{lemma}
    \begin{lemma}
  \label{lem1bis}
  Consider $P=\left[ \begin{smallmatrix} P_{11} & P_{12} \\[1mm] P_{21} & P_{22} \end{smallmatrix} \right]$ where  $P_{11}$ and $P_{22}$ are square and $P_{22}$ is non-singular. Then,
  \bea
  \label{det11}
  \det\,P=\det\, P_{22}\,\cdot\,\det(P_{11}-P_{12} P_{22}^{-1} P_{21}^\tra).
  \eea
  \end{lemma}
}
We now introduce the so-called Generalised Discrete Algebraic Riccati Equation GDARE($\Sigma$), defined as
\bea
\label{gdare}
X= A^\tra\,X\,A-(A^\tra\,X\,B+S)(R+B^\tra\,X\,B)^\dagger(B^\tra\,X\,A+S^\tra)+Q.
\eea

The algebraic equation (\ref{gdare}) subject to the constraint 
\bea
\label{kercond}
\ker (R+B^\tra\,X\,B) \subseteq \ker (A^\tra\,X\,B+S)
\eea
 is usually referred to as Constrained Generalised Discrete Algebraic Riccati Equation CGDARE($\Sigma$):
\bea
\label{cgdare}
\left\{ \begin{array}{ll}
X= A^\tra\,X\,A-(A^\tra\,X\,B+S)(R+B^\tra\,X\,B)^\dagger(B^\tra\,X\,A+S^\tra)+Q \\
\ker (R+B^\tra\,X\,B) \subseteq \ker (A^\tra\,X\,B+S) \end{array} \right.
\eea

 It is obvious that CGDARE($\Sigma$) constitutes a generalisation of the classic Discrete Riccati Algebraic Equation DARE($\Sigma$)
\bea
\label{dare}
X= A^\tra\,X\,A-(A^\tra\,X\,B+S)(R+B^\tra\,X\,B)^{-1}(B^\tra\,X\,A+S^\tra)+Q,
\eea
in the sense that any solution of DARE($\Sigma$) is also a solution of CGDARE($\Sigma$) but the {\em vice-versa} is not true in general. Importantly, however, the inertia of $R+B^\tra\,X\,B$ is independent of the particular solution of the CGDARE($\Sigma$), \cite[Theorem 2.4]{Stoorvogel-S-98}. This implies that a given CGDARE($\Sigma$) cannot have one solution $X=X^\tra$ such that $R+B^\tra X\,B$ is non-singular and another solution $Y=Y^\tra$ for which $R+B^\tra Y\,B$ is singular. As such, {\bf i)} if $X$ is a solution of DARE($\Sigma$), then all solutions of CGDARE($\Sigma$) will also satisfy DARE($\Sigma$) and, {\bf ii)} if $X$ is a solution of 
CGDARE($\Sigma$) such that $R+B^\tra\,X\,B$ is singular, then DARE($\Sigma$) does not admit solutions. \\

{To simplify the notation, for any $X=X^\tra\in  \real^{n \times n}$ we define
\beann
R_X & \defi  & R+B^\tra\,X\,B\\
S_X & \defi & A^\tra\,X\,B+S\\
K_X & \defi  &(R+B^\tra\,X\,B)^\dagger \, (B^\tra\,X\,A+S^\tra)=R_X^\dagger S_X^\tra\\
A_X&  \defi & A-B\,K_X
\eeann
%
%
%
%
so that (\ref{kercond}) can be written as $\ker R_X \subseteq \ker S_X$.%
}

\section{GDARE and the extended symplectic pencil}
{In this section we adapt the analysis carried out in \cite{Ferrante-W-07} for standard discrete algebraic Riccati equations to the case of CGDARE($\Sigma$).} 
Consider the so-called extended symplectic pencil $N-z\,M$, where
\beann
M\defi\left[ \begin{array}{ccc} 
I_n & O & O \\
O & -A^\tra & O \\ 
O & -B^\tra & O \end{array} \right] \qquad \textrm{and} \qquad
N\defi \left[ \begin{array}{ccc} 
A & O & B \\
Q & -I_n & S \\ 
S^\tra & O & R \end{array} \right].  
\eeann
{ This is an extension that may be reduced to the symplectic structure (see \cite{Wim,Ferrante-L-98}) when the matrix $R$ is invertible.}
We begin by giving a necessary and sufficient condition for $N$ to be singular. We will also show that, unlike the case in which the pencil $N-z\,M$ is regular, {the singularity of $N$ is not equivalent to the fact that the matrix pencil $N-z\,M$ has a generalised eigenvalue at zero.}

\begin{lemma}
\label{lem31}
Matrix $N$ is singular if and only if at least one of the two matrices $R$ and $A-B\,R^\dagger\,S^\tra$ is singular. 
\end{lemma}
\proof
First note that $N$ is singular if and only if such is $\bsmat A & B \\[1mm] S^\tra & R \esmat$.
To see this fact, consider the  left null-spaces.
Clearly, $\bmat{ccc} v_1\tp & v_2\tp & v_3\tp \emat\,N=0$, if and only if  $v_2=0$ and $\bmat{cc} v_1\tp & v_3\tp \emat\,\bsmat A & B \\[1mm] S^\tra & R \esmat=0$.\\ Now, if $R$ is singular, a non-zero vector $v_3$ exists such $v_3\tp\,R=0$. {Since from {\em (1)} in Lemma \ref{lem1} applied to the Popov matrix $\bsmat Q & S \\[1mm] S^\tra & R \esmat$ the subspace inclusion
$\ker R \subseteq \ker S$ holds}, we have also $\bmat{cc} 0\; & v_3\tp \emat \bsmat A & B \\[1mm] S^\tra & R \esmat=0$.
{If $R$ is invertible but $A-B\,R^\dagger\,S^\tra=A-B\,R^{-1}\,S^\tra$ is singular, {from (\ref{det11}) in Lemma \ref{lem1bis}} matrix $\bsmat A & B \\[1mm] S^\tra & R \esmat$ is singular, and therefore so is $N$.
{\em Vice-versa}, if both $R$ and $A-B\,R^{-1}\,S^\tra$ are non-singular, $\bsmat A & B \\[1mm] S^\tra & R \esmat$ 
is non-singular in view of {(\ref{det11}) in Lemma \ref{lem1bis}}. Thus, $N$ is invertible.\endproof}

 \ \\[-2mm]
 
The following theorem (see \cite{Ferrante-N-12-sub2} for a proof) presents a useful decomposition of the extended symplectic pencil that parallels the classic one -- see e.g. \cite{Ferrante-W-07} -- which is valid in the case in which the pencil $N-z\,M$ is regular.

\begin{theorem}
\label{the0}
Let $X$ be a symmetric solution of CGDARE($\Sigma$). Let also $K_X$ be the associated gain and $A_X$ be the associated closed-loop matrix. 
Two invertible matrices $U_X$ and $V_X$ of suitable sizes exist such that 
\bea
\label{decomp}
U_X\,(N-z\,M)\,V_X=\left[ \begin{array}{ccc}
A_X-z\,I_n & O & B \\
O & I_n-z\,A_X^\tra & O \\
O & -z\,B^\tra &  R_X \end{array} \right].
\eea
\end{theorem}

From Theorem \ref{the0} we find that if $X$ is a solution of CGDARE($\Sigma$), in view of the triangular structure obtained above we have
\bea
\label{det}
\det(N-z\,M)=(-1)^n \cdot \det(A_X-z\,I_n)\cdot \det(I_n-z\,A_X^\tra) \cdot \det R_X.
\eea

When $R_X$ is non-singular, the dynamics represented by this matrix pencil are decomposed into a part governed by the generalised eigenstructure of $A_X-z\,I_n$, a part governed by the finite generalised eigenstructure of $I_n-z\,A_X^\tra$, and a part which corresponds to the dynamics of the eigenvalues at infinity. 
When $X$ is a solution of DARE($\Sigma$), the generalised eigenvalues\footnote{Recall that a generalised eigenvalue of a matrix pencil $N-z\,M$ is a value of $z \in \complex$ for which the rank of the matrix pencil $N-z\,M$ is lower than its normal rank.}  of {$N\,z-M$} are given by the eigenvalues of $A_X$, the reciprocal of the non-zero eigenvalues of $A_X$, and {a generalised eigenvalue} at infinity whose algebraic  multiplicity is equal to $m$ plus the algebraic multiplicity of the eigenvalue of $A_X$ at the origin.
 The matrix pencil $I_n-z\,A_X^\tra$ has no generalised eigenvalues at $z=0$. This means that $z=0$ is a generalised eigenvalue of the matrix pencil $U_X\,(N-z\,M)\,V_X$ if and only if it is a generalised eigenvalue of the matrix pencil $A_X-z\,I_n$, because certainly 
 $z=0$ cannot cause the rank of $I_n-z\,A_X^\tra$ 
 to be smaller than its normal rank and because the normal rank of $N-z\,M$ is $2\,n+m$. This means that the Kronecker eigenstructure of the eigenvalue at the origin of $U_X\,(N-z\,M)\,V_X$ coincides with the Jordan eigenstructure of the eigenvalue at the origin of the closed-loop matrix $A_X$. Since the generalised eigenvalues of  $N-z\,M$ do not depend on the particular solution $X=X^\tra$ of CGDARE($\Sigma$), the same holds for the generalised eigenvalues and the Kronecker structure of $U_X\,(N-z\,M)\,V_X$ for any non-singular $U_X$ and $V_X$. Therefore, the nilpotent structure of the closed-loop matrix $A_X$ -- which is the Jordan eigenstructure of the generalised eigenvalue at the origin of $A_X$ -- if any, is independent of the particular solution $X=X^\tra$ of CGDARE($\Sigma$). Moreover, since
 \bea
 \label{newN}
U_X\,N\,V_X=\left[ \begin{array}{ccc}
A_X & O & B \\
O & I_n & O \\
O & O &  R_X \end{array} \right],
\eea
we see that, when $R_X$ is invertible, $N$ is singular if and only if $A_X$ is singular. { Since from Lemma \ref{lem31} matrix
$N$ is singular if and only if at least one of the two matrices $R$ and $A-B\,R^\dagger\,S^\tra$ is singular, we also have the following result. }

\begin{lemma}{{\bf (see e.g. \cite{Ferrante-04})}}
\label{prep1}
Let $R_X$ be invertible. Then, $A_X$ is singular if and only if at least one of the two matrices $R$ and $A-B\,R^\dagger\,S^\tra$ is singular.
\end{lemma}

However, when the matrix $R_X$ is singular, it is no longer true that 
 $A_X$ is singular if and only if $R$ or $A-B\,R^\dagger\,S^\tra$ is singular. Indeed, (\ref{newN}) shows that the algebraic multiplicity of the eigenvalue at the origin of $N$ is equal to the sum of the algebraic multiplicities of the eigenvalue at the origin of $A_X$ and $R_X$. Therefore, the fact that $N$ is singular does not necessarily imply that $A_X$ is singular. 
 {Indeed, Lemma \ref{prep1} can be generalised to the case where $R_X$ is possibly singular as follows.}

\begin{proposition}
\label{prep2}
The closed-loop matrix $A_X$ is singular if and only if $\rank R < \rank R_X$ or $A-B\,R^\dagger\,S^\tra$ is singular.
\end{proposition}
\proof Given a square matrix $Z$, let us denote by $\mu(Z)$ the algebraic multiplicity of its eigenvalue at the origin. Then, we know from (\ref{newN}) that $\mu(N)=\mu \left( \bsmat A & B \\[1mm] S^\tra & R \esmat \right)=\mu(A_X)+\mu(R_X)$. Consider a basis in the input space that isolates the invertible part of $R$. In other words, in this basis $R$ is written as $R=\bsmat R_1 & O \\[1mm] O & O \esmat$ where $R_1$ is invertible, while $B=\bmat{cc} B_1 & B_2 \emat$ and $S=\bmat{cc} S_1 & O \emat$ are partitioned accordingly. It follows that $\mu\left( \bsmat A & B \\[1mm] S^\tra & R \esmat \right)=\mu(R)+\mu\left( \bsmat A & B_1 \\[1mm] S_1^\tra & R_1 \esmat \right)$. As such,
\bea
\label{mu}
\mu(A_X)=\mu \left( \bmat{cc} A & B \\ S^\tra & R \emat \right)-\mu(R_X)=\mu \left( \bmat{cc} A & B_1 \\ S_1^\tra & R_1 \emat \right)+\mu(R)-\mu(R_X).
\eea
First, we show that if $\rank R < \rank R_X$, then $A_X$ is singular.  Since  $\rank R < \rank R_X$, then obviously $\mu(R)>\mu(R_X)$, so that (\ref{mu}) gives $\mu(A_X)>0$. \\ 
Let now $A-B\,R^\dagger\,S^\tra$ be singular, and let $\rank R = \rank R_X$. From (\ref{mu}) we find that $\mu(A_X)=\mu \left( \bsmat A & B_1 \\[1mm] S_1^\tra & R_1 \esmat\right)$. However, $A-B\,R^\dagger\,S^\tra=A-B_1\,R_1^{-1}\,S_1^\tra$. If $A-B\,R^\dagger\,S^\tra$ is singular, there exists a non-zero vector $k$ such that $\bmat{cc} k^\tra & -k^\tra\,B_1\,R_1^{-1}\emat \bsmat A & B_1 \\[1mm] S_1^\tra & R_1 \esmat=0$. Hence, $\mu \left(\bsmat A & B_1 \\[1mm] S_1^\tra & R_1 \esmat\right)>0$, and therefore also $\mu(A_X)>0$.  \\
To prove that the converse is true, it suffices to show that if $A-B\,R^\dagger\,S^\tra$ is non-singular
and $\rank R = \rank R_X$, then $A_X$ is non-singular. To this end, we observe that $\rank R=\rank R_X$ is equivalent to $\mu(R)=\mu(R_X)$ because $R$ and $R_X$ are symmetric. Thus, in view of (\ref{mu}), it suffices to show that if $A-B\,R^\dagger\,S^\tra$ is non-singular, then $\mu \left(\bsmat A & B_1 \\[1mm] S_1^\tra & R_1 \esmat\right)=0$. Indeed, {assume that  $A-B\,R^\dagger\,S^\tra=A-B_1\,R_1^{-1}\,S_1^\tra$ is non-singular, and} take a vector $\bsmat v_1^\tra & v_2^\tra \esmat$ such that $\bsmat v_1^\tra & v_2^\tra \esmat\bsmat A & B_1 \\[1mm] S_1^\tra & R_1 \esmat =0$. Then, since $R_1$ is invertible we get
$v_2^\tra=-v_1^\tra\,B_1\,R_1^{-1}$ and $v_1^\tra\,(A-B_1\,R_1^{-1}\,S_1^\tra)=0$. Hence, $v_1=0$ since $A-B_1\,R_1^{-1}\,S_1^\tra$ is non-singular, and therefore also $v_2=0$.
\endproof

\begin{remark}
 {\em
 We recall that $\mu(R_X)$ is invariant for any {symmetric} solution $X$ of CGDARE($\Sigma$), \cite{Stoorvogel-S-98}.
 Hence, as a direct consequence of (\ref{mu}), we have that $\mu(A_X)$ is the same for any {symmetric} solution $X$ of CGDARE($\Sigma$).  This means, in particular,  that  the closed-loop matrix corresponding to a given {symmetric} solution of CGDARE($\Sigma$) is singular if and only if the closed-loop matrix corresponding to any other {symmetric} solution of CGDARE($\Sigma$) is singular.
In the next section we show that a stronger result holds: when present, the zero eigenvalue has the same Jordan structure for any pair $A_X$ and $A_Y$ of closed-loop matrices corresponding to any pair $X,Y$ of {symmetric} solutions of CGDARE($\Sigma$). Moreover, the generalised eigenspaces corresponding to the zero eigenvalue of $A_X$ and $A_Y$ coincide. The restriction of $A_X$ and $A_Y$ to this generalised eigenspace also coincide. Finally, $X$ and $Y$ coincide along this generalised eigenspace. 
}
\end{remark}

\section{The subspace where all solutions coincide}
Given a solution $X=X^\tra$ of CGDARE($\Sigma$), we denote by $\gU$ the generalised eigenspace corresponding to the eigenvalue at the origin of $A_X$, i.e., $\gU \defi \ker (A_X)^n$. {Notice that, in principle, $\gU$ could depend on the particular solution $X$.
{In this section, and in particular in Theorem \ref{main},} we want to prove not only that $\gU$ does {\em not} depend on the particular solution $X$, but also that
 all solutions of CGDARE($\Sigma$) are coincident along $\gU$. In other words, given two solutions 
 $X=X^\tra$ and $Y=Y^\tra$ of CGDARE($\Sigma$), we show that $\ker (A_X)^n=\ker (A_Y)^n$ and, given a basis matrix\footnote{Given a subspace $\gS$, a basis matrix $S$ of $\gS$ is such that $\ima S=\gS$ and $\ker S=\{0\}$.} $U$ of the subspace $\gU=\ker (A_X)^n=\ker (A_Y)^n$,} the change of coordinate matrix $T=[\,U\;\;\;U_c\,]$ yields
 \bea
 T^{-1}\,X\,T=\left[ \begin{array}{cc} X_{11} & X_{12} \\ X_{12}^\tra & X_{22} \end{array} \right] \quad \textrm{and} \quad 
T^{-1}\,Y\,T=\left[ \begin{array}{cc} X_{11} & X_{12} \\ X_{12}^\tra & Y_{22} \end{array} \right]. \label{alph}
\eea

{We begin by presenting a first simple result.
\begin{lemma}
Two symmetric solutions  $X$ and $Y$ of CGDARE($\Sigma$) are coincident along the subspace $\gU$ if and only if 
$\gU \subseteq \ker (X-Y)$. 
\end{lemma}
\proof 
Suppose $X$ and $Y$ are coincident along the subspace $\gU$, and are already written in the basis defined by $T$ in (\ref{alph}). In this basis $\gU$ can be written as $\gU=\ima \left[ \begin{smallmatrix} I \\[1mm] O \end{smallmatrix} \right]$. If (\ref{alph}) holds, then we can write $X-Y=\bsmat O & O \\[1mm] O & \star \esmat$. Then, $(X-Y)\,\gU=\bsmat O & O \\[1mm] O & \star \esmat\left[ \begin{smallmatrix} I \\[1mm] O \end{smallmatrix} \right]=\{0\}$. {\em Vice-versa}, if $(X-Y)\,\gU=\{0\}$ and we write $X-Y=\bsmat \Delta_{11} & \Delta_{12} \\[1mm] \Delta_{12}^\tra & \Delta_{22} \esmat$, we find that $\bsmat \Delta_{11} & \Delta_{12} \\[1mm] \Delta_{12}^\tra & \Delta_{22} \esmat\left[ \begin{smallmatrix} I \\[1mm] O \end{smallmatrix} \right]=\{0\}$ implies $\Delta_{11}=0$ and $\Delta_{12}=0$.
\endproof
}
\ \\


{We now present two results that will be useful to prove Theorem \ref{main}.} Let $X=X^\tra\in \real^{n \times n}$. Similarly to \cite{Ferrante-W-07}, we define the function
 \bea
 \label{gdaredef}
{\cal D}(X)\defi X-A^\tra\,X\,A+(A^\tra\,X\,B+S)(R+B^\tra\,X\,B)^\dagger(B^\tra\,X\,A+S^\tra)-Q.
\eea
If in particular $X=X^\tra$ is a solution of GDARE($\Sigma$), then $\gD(X)=0$.
{Recall that  we have defined $R_X = R+B^\tra\,X\,B$, $S_X = A^\tra\,X\,B+S$ and $R_Y = R+B^\tra\,Y\,B$, $S_Y \defi A^\tra\,Y\,B+S$.}

 \begin{lemma}
 \label{41}
 Let $X=X^\tra\in \real^{n \times n}$ and $Y=Y^\tra\in \real^{n \times n}$ be such that (\ref{kercond}) holds, i.e., 
 \bea
 \ker R_X  \subseteq  \ker S_X \label{kercondX} \\
 \ker R_Y  \subseteq  \ker S_Y.\label{kercondY}
 \eea
 Let $A_X = A-B\,K_X$ with $K_X = R_X^\dagger\,S_X^\tra$ and  $A_Y = A-B\,K_Y$ with $K_Y = R_Y^\dagger\,S_Y^\tra$. 
Moreover, let us define the difference $\Delta\defi X-Y$. Then, 
\bea
\label{ionescu}
{\cal D}(X)-{\cal D}(Y)=\Delta-A_Y^\tra\,\Delta\,A_Y+A_Y^\tra\,\Delta\,B\,R_X^\dagger\, B^\tra\,\Delta\, A_Y.
\eea
\end{lemma}
The proof can be found in \cite[p.382]{Abou-Kandil-FIJ-03}.

  {The following lemma is the counterpart of Lemma 2.2 in \cite{Ferrante-W-07} where the standard DARE was considered.}  
  
\begin{lemma}
\label{lemWF}
 Let $X=X^\tra\in \real^{n \times n}$ and $Y=Y^\tra\in \real^{n \times n}$ be such that (\ref{kercondX}-\ref{kercondY}) hold.  {Let $\Delta=X-Y$.} Then,
\bea
{\cal D}(X)-{\cal D}(Y)=\Delta-A_Y^\tra\,\Delta \, A_X.
\eea
\end{lemma}
\proof
First, notice that
\beann
A_Y^\tra \Delta\, B = [A^\tra-(A^\tra Y \,B +S)\,R_Y^\dagger B^\tra] \Delta \,B.
\eeann
We now show that $\ker R_X \subseteq \ker (A_Y^\tra \Delta\, B)$. To this end, let $P_X$ be a basis of the null-space of $R_X$. Hence, $(R+B^\tra X B)P_X=0$.  Then,
\beann
A_Y^\tra\,\Delta\,B\,P_X  &= & \left(A^\tra-(A^\tra\,Y\,B+S)\,R_Y^\dagger\,B^\tra\right)\,(X-Y)\,B\,P_X \\
 &= & A^\tra\,X\,B\,P_X-(A^\tra\,Y\,B+S)\,R_Y^\dagger\,B^\tra\,X\,B\,P_X -A^\tra\,Y\,B\,P_X \\
&&  +(A^\tra Y \,B +S)\,R_Y^\dagger\,B^\tra\,Y\,B\,P_X \\
&&  +(A^\tra Y \,B +S)\,R_Y^\dagger\,R\,P_X- (A^\tra Y \,B +S)\,R_Y^\dagger\,R\,P_X \\
& =  &A^\tra\,X\,B\,P_X+(A^\tra Y \,B +S)\,R_Y^\dagger\,R_Y\,P_X-A^\tra\,Y\,B\,P_X \\
& = & A^\tra\,X\,B\,P_X+S_Y\,P_X-A^\tra\,Y\,B\,P_X=(A^\tra\,X\,B+S)\,P_X,
\eeann
which is zero since $\ker R_X \subseteq \ker S_X$ in view of { (\ref{kercondX}) in Lemma \ref{41}}. Now we want to prove that
\bea
\label{eqQ}
A_Y^\tra \Delta \,(A_Y-A_X)  =   A_Y^\tra \, \Delta\, B\,R_X^\dagger \, B^\tra\, \Delta\,A_Y. 
\eea
Consider the term
\bea
\label{ok}
A_Y^\tra \Delta (A_Y-A_X)  =  A_Y^\tra \Delta\, B \,(R_X^\dagger S_X-R_Y^\dagger S_Y).
\eea
{
Since $R_X^\dagger R_X$ is an orthogonal projection that projects onto $\ima R_X^\tra =\ima R_X$, we have $\ker R_X=\ima (I_m-R_X^\dagger R_X)$. Since as we have shown $\ker R_X \subseteq \ker (A_Y^\tra \Delta\, B)$, from$\ker R_X=\ima (I_m-R_X^\dagger R_X)$ we also have}
$A_Y^\tra \Delta\, B \,(I_m-R_X^\dagger R_X)=0$, which means that $A_Y^\tra \Delta\, B\,R_X^\dagger\, R_X=A_Y^\tra \,\Delta\,B$. We use this fact on (\ref{ok}) to get
\bea
A_Y^\tra \Delta (A_Y\!-\!A_X) & \!=\! & A_Y^\tra \Delta\, B\,R_X^\dagger [ (B^\tra X A+S)-R_X\,R_Y^\dagger (B^\tra Y A+S)] \nonumber \\
   & \!=\! &A_Y^\tra \Delta\, B\,R_X^\dagger [ (B^\tra X A\!+\!S\!-\!B^\tra\,Y\,A\!+\!B^\tra Y \,A)\!-\!R_X\,R_Y^\dagger (B^\tra Y A\!+\!S)]  \nonumber \\
   & \!=\! & A_Y^\tra \Delta\, B\,R_X^\dagger [ B^\tra \Delta\,A+(I_m-R_X\,R_Y^\dagger) (B^\tra Y A+S)]. \label{final1}
 \eea
Since $R_X=R+B^\tra X \,B-B^\tra Y \,B+B^\tra Y \,B=R_Y+B^\tra\,\Delta\,B$, eq. (\ref{final1}) becomes
\beann
A_Y^\tra \Delta (A_Y-A_X)  =  A_Y^\tra \Delta\, B\,R_X^\dagger [B^\tra \Delta\,A + (I_m-R_Y\,R_Y^\dagger-B^\tra \Delta\,B\,R_Y^\dagger)(B^\tra Y A+S)] \\
  =  A_Y^\tra \Delta B\,R_X^\dagger B^\tra \Delta\, (A - B\,R_Y^\dagger)(B^\tra Y A+S)= \Delta B\,R_X^\dagger B^\tra \Delta\,A_Y,
\eeann
since from Lemma \ref{lem1} $(I_m-R_Y\,R_Y^\dagger)(B^\tra Y A+S)=0$ from $\ker R_Y \subseteq \ker (A^\tra Y\,B+S)$. Eq. (\ref{eqQ}) follows by recalling that $A_Y=A - B\,R_Y^\dagger\,S_Y$. Plugging (\ref{eqQ}) into (\ref{ionescu}) yields
\beann
{\cal D}(X)-{\cal D}(Y)=\Delta-A_Y^\tra\,\Delta A_Y+A_Y^\tra \Delta (A_Y-A_X)=\Delta-A_Y^\tra\,\Delta A_X.
\eeann
\endproof

Now we are ready to prove the main result of this section. This result {extends the analysis of Proposition 2.1 in \cite{Ferrante-W-07} to solutions of CGDARE($\Sigma$).}

\begin{theorem}
\label{main}
Let $\gU = \ker (A_X)^n$ denote the generalised eigenspace corresponding to the eigenvalue at the origin of $A_X$. Then
\begin{enumerate}
\item All solutions of CGDARE($\Sigma$) are coincident along $\gU$, i.e., given two solutions $X$ and $Y$ of CGDARE($\Sigma$),
\[
(X-Y)\,\gU=\{0\};
\]
\item $\gU$ does not depend on the solution $X$ of CGDARE($\Sigma$), i.e., given two solutions $X$ and $Y$ of CGDARE($\Sigma$), there holds
\beann
\ker (A_X)^n=\ker (A_Y)^n.
\eeann
\end{enumerate}
\end{theorem}
\proof Let us prove {\em (1)}. Consider a non-singular $T \in \real^{n \times n}$. Define the new quintuple
\beann
\tilde{A} \defi T^{-1}\,A\,T, \qquad \tilde{B}\defi T^{-1}\,B, \quad \tilde{Q}\defi T^\tra\,Q\,T, \quad  \tilde{S}\defi T^\tra S, \quad  \tilde{R}\defi R.
\eeann
It is straightforward to see that $X$ satisfies GDARE($\Sigma$) with respect to $(A,B,Q,R,S)$ if and only if $\tilde{X}\defi T^\tra X\,T$ satisfies GDARE($\Sigma$) with respect to $(\tilde{A},\tilde{B},\tilde{Q},\tilde{R},\tilde{S})$, which for the sake of simplicity is denoted by $\tilde{\gD}$, so that $\tilde{\gD}(\tilde{X})=0$. The closed-loop matrix in the new basis is related to the closed-loop matrix in the original basis by
{ 
$$
\tilde{A}_{\tilde{X}} =  \tilde{A}-\tilde{B}\,(\tilde{R}+\tilde{B}^\tra \tilde{X} \,\tilde{B})^\dagger (\tilde{B}^\tra \tilde{X} \,\tilde{A}+\tilde{S}^\tra) =T^{-1}\,A_X\,T.
$$
}
Moreover, if $\tilde{\gU}=\ker (\tilde{A}_{\tilde{X}})^n$, then $\tilde{\gU}=T^{-1}\, \gU$ since $(\tilde{A}_{\tilde{X}})^n \tilde{\gU}=0$ is equivalent to $T^{-1} (A_X)^n T\,\tilde{\gU}=T^{-1} (A_X)^n\,\gU=0$.
We choose an orthogonal change of coordinate matrix $T$ as $T=[\,U\;\;\;U_c\,]$, where $U$ is a basis matrix of $\gU$. In this new basis
\beann
\tilde{A}_{\tilde{X}}=T^{-1}\,A_X\,T = \left[ \begin{array}{cc} U & U_c \end{array} \right]^\tra A_X \left[ \begin{array}{cc} U & U_c \end{array} \right] \\
 =  \left[ \begin{array}{cc} U^\tra A_X\,U & \star \\ U_c^\tra A_X\,U & \star  \end{array} \right]=
\left[ \begin{array}{cc} U^\tra A_X\,U & \star \\ O & U_c^\tra A_X\,U_c  \end{array} \right],
\eeann
where the zero in the bottom left corner is due to the fact that the rows of $U_c^\tra A_X$ are orthogonal to the columns of $U$. Moreover, the submatrix $N_0 \defi U^\tra A_X\,U$ is nilpotent with the same nilpotency index\footnote{With a slight abuse of nomenclature, we use the term {\em nilpotency index} of a matrix $M$ to refer to the smallest integer $\nu$ for which $\ker (M)^\nu=\ker (M)^{\nu+1}$, which is defined also when $M$ is not nilpotent.} of $A_X$. 
Notice also that $H_X \defi U_c^\tra A_X\,U_c$ is non-singular.
Let $\tilde{X}$ be a solution of CGDARE($\tilde{\Sigma}$) in this new basis, and let it be partitioned as
\beann
\tilde{X}=\left[ \begin{array}{cc} \tilde{X}_{11} & \tilde{X}_{12} \\ \tilde{X}_{12}^\tra & \tilde{X}_{22} \end{array} \right],
\eeann
where $\tilde{X}_{11}$ is $\nu \times \nu$, with $\nu=\textrm{dim} \,\gU$. Consider another solution $\tilde{Y}$ of CGDARE($\tilde{\Sigma}$), partitioned as $Y=\bsmat \tilde{Y}_{11} & \tilde{Y}_{12} \\[1mm] \tilde{Y}_{12}^\tra & \tilde{Y}_{22} \esmat$. Let $\Delta \defi \tilde{X}-\tilde{Y}$ be partitioned in the same way. Since $\tilde{X}$ and $\tilde{Y}$ are both solutions of CGDARE($\tilde{\Sigma}$), we get $\tilde{\gD}(\tilde{X})=\tilde{\gD}(\tilde{Y})=0$. Thus, in view of Lemma \ref{lemWF}, there holds
\bea
\label{eqalpha}
\Delta-\tilde{A}_{\tilde{Y}}^\tra \,\Delta \,\tilde{A}_{\tilde{X}}=0.
\eea
If $\Delta$ is partitioned as $\Delta=[\,\Delta_1\;\;\;\Delta_2\,]$ where $\Delta_1$ has $\nu$ columns, eq. (\ref{eqalpha}) becomes
\[
\left[ \begin{array}{cc} \Delta_1 & \Delta_2 \end{array} \right]-\tilde{A}_{\tilde{Y}}^\tra \left[ \begin{array}{cc} \Delta_1 & \Delta_2 \end{array} \right]\left[ \begin{array}{cc} N_0 & \star \\ O & H_X  \end{array} \right]=
\left[ \begin{array}{cc} \Delta_1 -\tilde{A}_{\tilde{Y}}^\tra \Delta_1\,N_0 & \star \end{array} \right]=0,
\]
from which we get $\Delta_1=\tilde{A}_{\tilde{Y}}^\tra\, \Delta_1\,N_0$. Thus,
\beann
\Delta_1 = \tilde{A}_{\tilde{Y}}^\tra \Delta_1\,N_0=(\tilde{A}_{\tilde{Y}}^\tra)^2 \Delta_1\,N_0^2= \ldots =
(\tilde{A}_{\tilde{Y}}^\tra)^n \Delta_1\,(N_0)^n,
\eeann
which is equal to zero since $(N_0)^n$ is the zero matrix. Hence, $\Delta_1=0$. Thus, we have also
\[
\Delta\,{\gU}=\bmat{cc} O \; &\; \star \emat \left( \ima\bmat{c} I \\ O \emat\right)=\{0\}.
\]
Since $\Delta$ is symmetric, we get
\beann
\tilde{X}-\tilde{Y}=\left[ \begin{array}{cc} \tilde{X}_{11} & \tilde{X}_{12} \\ \tilde{X}_{12}^\tra & \tilde{X}_{22} \end{array} \right]-\left[ \begin{array}{cc} \tilde{Y}_{11} & \tilde{Y}_{12} \\ \tilde{Y}_{12}^\tra & \tilde{Y}_{22} \end{array} \right]=\left[ \begin{array}{cc} O & O \\ O & \tilde{X}_{22}-\tilde{Y}_{22} \end{array} \right],
\eeann
which leads to $\tilde{X}_{11}=\tilde{Y}_{11}$ and $\tilde{X}_{12}=\tilde{Y}_{12}$. \\

Let us prove {\em (2)}. 
Since $\ker R_Y$ coincides with $\ker R_X$ by virtue of \cite[Theorem 4.3]{Ferrante-N-12-sub}, we find
\bea
A_X-A_Y & =&  B\,(R_Y^\dagger S_Y^\tra-R_X^\dagger \,S_X^\tra) \nonumber \\
& =  & B\,R_Y^\dagger (S_Y^\tra-R_Y\,R_X^\dagger \,S_X^\tra). \label{sec}
\eea  
Plugging 
\bea
\label{eqsy}
S_Y^\tra=B^\tra \,Y\,A+S^\tra=B^\tra \,\Delta\,A+S^\tra+B^\tra\,X\,A =B^\tra\,\Delta\,A+S_X^\tra
\eea
and
\bea
\label{eqry}
R_Y =R+ B^\tra \,Y\,B-B^\tra\,X\,B+B^\tra\,X\,B=R_X+B^\tra \,\Delta\,B
\eea
into (\ref{sec}) yields
\beann
A_X-A_Y & = & B\,R_Y^\dagger (B^\tra\,\Delta\,A-B^\tra\,\Delta\,B\,R_X^\dagger \,S_X^\tra) \\
& = & B\,R_Y^\dagger B^\tra\,\Delta\,A_X.
\eeann
This means that the identity
\[
A_X-A_Y=B\,R_Y^\dagger B^\tra\,\Delta\,A_X
\]
holds. By partitioning $\Delta=\bsmat O & \star \\[1mm]
O & \star \esmat$, we find that also $B\,R_Y^\dagger B^\tra\,\Delta=\bsmat O & \star \\[1mm]
O & \star \esmat$, so that
\beann
A_Y  &=&  A_X-B\,R_Y^\dagger B^\tra\,\Delta\,A_X \\
& = & \bmat{cc} N_0 & \star \\ O & H_X \emat -\bmat{cc} O \;&\; \star \\ O & \star \emat\bmat{cc} N_0 & \star \\ O & H_X \emat =\bmat{cc} N_0 & \star \\ O & H_Y \emat.
\eeann
Thus, $\ker (A_Y)^n \supseteq \ker (A_X)^n$. If we interchange the role of $X$ and $Y$, we obtain the opposite inclusion
$\ker (A_Y)^n \subseteq \ker (A_X)^n$. Notice, in passing, that this also implies that $H_Y$ is non-singular.
 \endproof

\section{The Generalised Riccati Difference Equation}
Consider the GRDE($\Sigma$) along with the terminal condition $X_T=P=P^\tra\ge 0$. Let us define
\beann
\gR(X) \defi A^\tra\,X\,A-(A^\tra\,X\,B+S)(R+B^\tra X\,B)^\dagger(B^\tra\,X\,A+S^\tra)+Q.
\eeann
With this definition, GRDE($\Sigma$) can be written as $X_{t}=\gR(X_{t+1})$. Moreover, GDARE($\Sigma$) can be written as
\beann
\gD(X)=X-\gR(X)=0.
\eeann

We have the following important result.

\begin{theorem}
\label{th51}
Let $X_\circ=X_\circ^\tra$ be a solution of CGDARE($\Sigma$). Let $\nu$ be the index of nilpotency of $A_{X_\circ}$. Moreover, let $X_t$ be a solution of (\ref{grde}-\ref{term}) and define $\Delta_t \defi X_t-X_\circ$. Then, for $\tau \ge \nu$, we have $\Delta_{T-\tau}\,\gU=\{0\}$.
\end{theorem}
\proof
Since $X_\circ=X_\circ^\tra$ is a solution of CGDARE($\Sigma$), we have $\gD(X_\circ)=0$. This is equivalent to saying that $X_\circ=\gR(X_\circ)$. From the definition of $\Delta_t$ we get in particular $\Delta_T=X_T-X_\circ$. With these definitions in mind, we find
\bea
\Delta_t  &=&  \gR(X_{t+1})-\gR(X_\circ)=X_{t+1}-\gD(X_{t+1})-X_\circ \nonumber \\
 &= & \Delta_{t+1}-\gD(X_{t+1})=\Delta_{t+1}-\gD(X_{t+1})+\gD(X_\circ) \nonumber  \\
 &=&  \Delta_{t+1}-[\gD(X_{t+1})-\gD(X_\circ)]. \label{pippo} 
\eea
However, we know from (\ref{ionescu}) that
\bea
\label{ionescu1}
&&{\cal D}(X_{t+1})-{\cal D}(X_\circ) \nonumber \\
&& \hspace{.3cm} = \Delta_{t+1}-A_{X_\circ}^\tra\,[\Delta_{t+1}-\Delta_{t+1}\,B\, (R+B^\tra X_{t+1} B)^\dagger B^\tra\,\Delta_{t+1} ]A_{X_\circ},
\eea
which, once plugged into (\ref{pippo}), gives
\bea
\Delta_t  &=&  \Delta_{t+1}-\Delta_{t+1}+A_{X_\circ}^\tra\,[\Delta_{t+1}+\Delta_{t+1}\,B\, (R+B^\tra X_{t+1} B)^\dagger B^\tra\,\Delta_{t+1} ]A_{X_\circ} \nonumber \\
 &=&  A_{X_\circ}^\tra\,[I_n-\Delta_{t+1}\,B\, (R+B^\tra X_{t+1} B)^\dagger B^\tra\,] \Delta_{t+1} A_{X_\circ}=F_{t+1}\,\Delta_{t+1} \,A_{X_\circ}, \label{Riccati}
\eea
where
\[
F_{t+1} \defi A_{X_\circ}^\tra-A_{X_\circ}^\tra\Delta_{t+1}\,B\, (R+B^\tra X_{t+1} B)^\dagger B^\tra.
\]
It follows that we can write
\bea
\Delta_{T-1}  &=&  F_T\,\Delta_T\,A_{X_\circ}, \nonumber \\
\Delta_{T-2} & =&  F_{T-1}\,\Delta_{T-1}\,A_{X_\circ}=F_{T-1}\,F_T\,\Delta_T\,(A_{X_\circ})^2,  \nonumber \\
 & \vdots &   \\
\Delta_{T-\tau} & = & \left(\prod_{i=T-\tau+1}^T F_i\right)\,\Delta_T\,(A_{X_\circ})^{\tau}. \label{ultima}
\eea 
This shows that for $\tau \ge \nu$ we have $\ker \Delta_{T-\tau} \supseteq \ker (A_{X_\circ})^n$. \endproof

\ \\[-2mm]

Now we show that the result given in Theorem \ref{th51} can be used to obtain a reduction for the generalised discrete-time Riccati difference equation. Consider the same basis induced by the change of coordinates used in Theorem \ref{main}, so that the first $\nu$ components of this basis span the subspace $\gU=\ker (A_X)^n$. The closed-loop matrix in this basis can be written as
\beann
{A}_{{X_\circ}} = \left[ \begin{array}{cc} N_0 \;& \star \\ O & Z \end{array} \right],
\eeann
where $N_0$ is nilpotent and $Z$ is non-singular. Hence, $({A}_{{X_\circ}})^{\nu}=\bsmat O & \star \\[1mm] O & Z^{\nu} \esmat$, where we recall that $\nu$ is the nilpotency index of $A_{X_\circ}$. By writing (\ref{ultima}) in this basis, for $\tau \ge \nu$ we find 
\beann
{\Delta}_{T-\tau}=\bmat{cc} \star \;&\; \star \\ \star \;&\; \star \emat \left[ \begin{array}{cc} O\; &\; \star \\ O \;&\; Z^{\tau} \end{array} \right]
=
\left[ \begin{array}{cc} O \;&\; \star \\ O \;&\; \star \end{array} \right]=\left[ \begin{array}{cc} O \;& \;O \\ O\; &\; \star \end{array} \right],
\eeann
where the last equality follows from the fact that ${\Delta}_{T-\tau}$ is symmetric. 

Now, let us rewrite the Riccati difference equation (\ref{Riccati}) as
\bea
\Delta_t=  A_{X_\circ}^\tra\,\Delta_{t+1} A_{X_\circ}     -A_{X_\circ}^\tra\,\Delta_{t+1}\,B (R+B^\tra X_{t+1} B)^\dagger B^\tra\, \Delta_{t+1} A_{X_\circ}.
\eea
For $t \le T-\nu$, we get $\Delta_t=\bsmat O & O \\[1mm] O & \Psi_t \esmat$, and the previous equation becomes
\beann
\bmat{cc} \! O \; &   O  \! \\  \! O \;  &   \Psi_t  \! \emat & = &  
\left[ \begin{array}{cc}  \! N_0^\tra   &   O \!  \\  \! \star  \! & \!  Z^\tra  \! \end{array} \right]\bmat{cc}  \! O \; &   O \!  \\ \!  O \; &   \Psi_{t+1}  \! \emat\left[ \begin{array}{cc} \!  N_0   & \; \star  \! \\  \! O   & \;  Z  \! \end{array} \right]\\
& & -
\left[ \begin{array}{cc}  \! N_0^\tra  &   O  \! \\ \!  \star  \! &  \! Z^\tra \!  \end{array} \right]\bmat{cc}  \! O  &   O  \! \\  \! O  & \Psi_{t+1}  \! \emat B\,(R+B^\tra X_{t+1}\,B)^\dagger B^\tra \bmat{cc}  \! O  &   O  \! \\ \!  O & \Psi_{t+1} \!  \emat\left[ \begin{array}{cc}  \! N_0  &  \; \star  \! \\ \!  O \!  &  \! Z  \! \emat \\
& = & 
\bmat{cc} \! O   &   O \!  \\  \! O  &  Z^\tra \!  \, \Psi_{t+1} \, Z \! \emat  \\
& & -
 \bmat{cc} \!  O   &   O  \! \\ \!  O   &   Z^\tra \, \Psi_{t+1} \!  \emat  \! \! \bmat{c} \!  B_1 \!  \\ \!  B_2 \!  \emat  \!  \! \left(R \! + \! \bmat{cc}  \! B_1^\tra & B_2^\tra \!  \emat  \! (\Delta_{t+1} \! + \! X_\circ) \! \bmat{c}  \! B_1  \! \\  \! B_2 \!  \emat \right)^\dagger\!  \bmat{cc}  \! B_1^\tra & B_2^\tra \!  \emat  \! \!   \bmat{cc} \!  O  &  O  \! \\  \! O  &  \Psi_{t+1} Z \! \emat\! \! .
\eeann
By partitioning $X_\circ$ as { $X_\circ=\bsmat X_{\circ,11} & X_{\circ,12} \\[1mm] X_{\circ,12}^\tra & X_{\circ,22} \esmat$}, we get
\beann
\bmat{cc} O & O \\ O & \Psi_t \emat  &\!=\!&   
\bmat{cc} O & O \\ O & Z^\tra \, \Psi_{t+1} \, Z\emat \!-\! 
 \bmat{cc} O & O \\ O & Z^\tra \, \Psi_{t+1} \emat \!\!\bmat{cc} \star & \star  \\ \star & B_2\,(R_0\!+\!B_2^\tra\,\Psi_{t+1}\,B_2)^\dagger \,B_2^\tra \emat\!\! \bmat{cc} O & O \\ O & \Psi_{t+1} \, Z\emat \\
 &\!=\!&
\bmat{cc} O & O \\ O & Z^\tra \, \Psi_{t+1} \, Z\emat -
 \bmat{cc} O & O \\ O & Z^\tra \, \Psi_{t+1}\,B_2\,(R_0+B_2^\tra\,\Psi_{t+1}\,B_2)^\dagger \,B_2^\tra\,\Psi_{t+1} \, Z \emat,
 \eeann 
 where $R_0 \defi R+B_2^\tra \,X_{\circ,22}\,B_2$.  Therefore, $\Psi_t$ satisfies the reduced homogeneous Riccati difference equation
 \bea
 \label{reduced}
 \Psi_t=Z^\tra \, \Psi_{t+1} \, Z- Z^\tra \, \Psi_{t+1}\,B_2\,(R_0+B_2^\tra\,\Psi_{t+1}\,B_2)^\dagger \,B_2^\tra\,\Psi_{t+1} \, Z.
 \eea
The associated generalised discrete Riccati algebraic equation is
\bea
\label{homog}
\Psi- Z^\tra \, \Psi \, Z+Z^\tra \, \Psi\,B_2\,(R_0+B_2^\tra\,\Psi\,B_2)^\dagger \,B_2^\tra\,\Psi \, Z=0.
 \eea
 Being homogeneous, this equation admits the solution $\Psi=0$. This fact has two important consequences:
 \begin{itemize}
 \item The closed-loop matrix associated with this solution is clearly $Z$, which is non-singular. On the other hand, we know that the nilpotent part of the closed-loop matrix is independent of the particular solution of CGDARE($\Sigma$) considered. This means that all solutions of (\ref{homog}) have a closed-loop matrix that is non-singular;
 \item Given a solution $\Psi$ of (\ref{homog}), the null-space of $R_0+B_2^\tra\,\Psi\,B_2$ coincides with the null-space of $R_0$, since the null-space of $R_0+B_2^\tra\,\Psi\,B_2$ does not depend on the particular solution of (\ref{homog}) and we know that the zero matrix is a solution of (\ref{homog}).
 \end{itemize}

 As a result of this discussion, it turns out that given a reference solution $X_\circ$ of CGDARE($\Sigma$), the solution of GDRE($\Sigma$) with terminal condition $X_T=P$ can be computed backward as follows:
 \begin{enumerate}
 \item For the first $\nu$ steps, i.e., from $t=T$ to $t=T-\nu$, $X_t$ is computed by iterating the GDRE($\Sigma$) starting from the terminal condition $X_T=P$;
 \item In the basis that isolates the nilpotent part of $A_X$, we have 
 \[
 \Delta_{T-\nu}=\bmat{cc} O & O \\ O & \Psi_{T-\nu} \emat.
 \]
 From $t=T-\nu-1$ to $t=0$, the solution of GDRE($\Sigma$) can be found iterating the reduced order GDRE in (\ref{reduced}) starting from the terminal condition $\Psi_{T-\nu}$.
 \end{enumerate}
 
 \begin{remark}
 {\em The advantage of using the reduced-order generalised difference Riccati algebraic equation (\ref{reduced}) consists in the fact that the closed-loop matrix of any solution of the associated generalised discrete Riccati algebraic equation is non-singular. Hence, when the reduced-order pencil given by the Popov triple $\left(Z,B_2,\bsmat 0 & 0 \\[1mm] 0 & R_0 \esmat\right)$ is regular, the solution of the 
 reduced-order generalised difference Riccati algebraic equation (\ref{reduced}) can also be computed in closed-form, using the results in \cite{Ferrante-N-06}. Indeed, consider a solution $\Psi$ of (\ref{homog}) with its non-singular closed-loop matrix $A_\Psi$ and let $Y$ be the corresponding solution of the closed-loop Hermitian Stein equation
 \bea
\label{stein}
A_\Psi\,Y\,A_\Psi^{\tra}-Y+B_2\,(R_0+B_2^{\tra}\,\Psi\,B_2)^{-1}B_2^{\tra}=0. 
\eea
 
 The set of solutions of the extended symplectic difference equation for the reduced system is parameterised in terms of $K_1,K_2 \in \mathbb{R}^{(n-\nu) \times (n-\nu)}$ as
\bea
\label{param}
\bmat{c} \Xi_{t} \\ \Lambda_t \\ \Omega_t \emat=
\left[ \begin{array}{c} \! I_{n-\nu} \! \\ \! \Psi \! \\ \! -K_\Psi \! \end{array} \right]  (A_\Psi)^t \,K_1 \!+\! \left[ \begin{array}{c} \! Y\,A_\Psi^{\tra} \! \\ \! (\Psi\,Y-I_{n-\nu})A_\Psi^{\tra} \! \\ \! -K_{\star} \! \end{array} \right] (A_\Psi^{\tra})^{T-t-1}\, K_2, & \!\quad  0 \le t \le T,\quad
\eea
where $K_\star \defi K_{\Psi}\,Y\,A_{\Psi}^\tra - (R_0+B_2^\tra \,\Psi\,B_2)^{-1}\,B_2^\tra$. 
The values of the parameter matrices $K_1$ and $K_2$ can be computed so that the terminal condition satisfies $X_T=I_n$ and $\Lambda_T=\Psi_{T-\nu}$. Such values exist because $A_\Psi$ is non-singular, and are given by
 \beann
K_1  &=&  (A_\Psi)^{-T}\left(I_{n-\nu}-Y\,(\Psi-\Psi_{T-\nu}) \right) \\
K_2 & = & \Psi-\Psi_{T-\nu}.
\eeann
 Then, the solution of (\ref{reduced}) is given by $\Psi_t=\Lambda_t\,\Xi_t^{-1}$. 
 }
 \end{remark}

 {
 \section{Concluding remarks}
 In this paper we have considered the generalised Riccati difference equation with a terminal condition which arises in finite-horizon LQ optimal control. We have shown in particular that it is possible to identify and deflate the singular part of such equation using the corresponding generalised algebraic Riccati equation. The two advantages of this technique are the reduction of the dimension of the Riccati equation at hand as well as the fact that the reduced problem is non-singular, and can therefore be handled with the standard tools of the finite-horizon LQ theory. }


\end{document}